\documentclass[preprint,12pt]{elsarticle}

\usepackage{amssymb}
\usepackage{amsmath}
\usepackage{amsthm}

\usepackage{color}

\usepackage{url}
\usepackage{breakurl}
\usepackage[breaklinks,hyperfootnotes=false]{hyperref}

\newcommand{\todo}[1]{}
\renewcommand{\todo}[1]{{\color{red}TODO: {{#1}}}}

\newcommand{\CA}{\mathsf{CA}}
\newcommand{\OA}{\mathsf{OA}}
\newcommand{\CAN}{\mathsf{CAN}}
\newcommand{\CAK}{\mathsf{CAK}}

\newcommand*{\helvetize}{\fontfamily{phv}\selectfont} 
\newcommand{\hely}{\textit{\helvetize y}}
\renewcommand{\beta}{\hely}
\newcommand{\bbeta}{\textbf{\hely}}
\renewcommand{\alpha}{\lambda}
\newcommand{\balpha}{\boldsymbol{\lambda}}

\newtheorem{theorem}{Theorem}
\newtheorem{corollary}{Corollary}
\newtheorem{definition}{Definition}

\newcommand\blfootnote[1]{%
  \begingroup
  \renewcommand\thefootnote{}\footnote{#1}%
  \addtocounter{footnote}{-1}%
  \endgroup
}

\begin{document}

\begin{frontmatter}



\title{
Determining Covering Array Numbers via Balanced Covering Arrays
} 


\author[SBA]{Irene Hiess}
\ead{ihiess@sba-research.org}
\author[SBA,AAAM]{Ludwig Kampel\corref{LK}} 
\cortext[LK]{Corresponding author}
\ead{lkampel@aaam.info}
\affiliation[SBA]{organization={SBA Research}, 
            addressline={Floragasse 7}, 
            city={Vienna},
            postcode={1040}, 
            country={Austria}}
\affiliation[AAAM]{organization={AAAM Research},
            city={Vienna},
            country={Austria}}

\begin{abstract}

In this article we determine five previously unknown covering array numbers (CANs).
We do so using properties of so called balanced covering arrays together with a computational result for these.
The balance properties allow us to generalize the (computational) non-existence result for balanced covering arrays to covering arrays.
Covering arrays are combinatorial designs that can be considered generalizations of orthogonal arrays,
when dropping the restriction that the considered $t$-tuples appear \emph{exactly} $\lambda$ times,
and instead require them to appear \emph{at least} $\lambda$ times.
While this generalization renders the existence of covering arrays trivial,
it raises the question for their optimality, respectively the smallest number of rows, the CAN, for which a certain covering array exists.
The CANs determined in this paper were tightly bound for decades, but remained ultimately unknown.

\end{abstract}


\end{frontmatter}

\section{Introduction}

Covering arrays (CAs) are combinatorial designs that can be introduced as generalizations of orthogonal arrays \cite{SloaneOA}.\blfootnote{\emph{Preprint submitted to Elsevier} \hfill \emph{October 20, 2025}}
A $\CA_{\lambda}(N;t,k,v)$ is an array with $N$ rows,
$k$ columns and entries from a $v$-ary alphabet, with the defining property that every $t$-tuple over the alphabet occurs at least $\lambda$ times in every $N\times t$ subarray.
The value $t$ is called the \emph{strength} of the CA, and $\lambda$ its index. When $\lambda = 1$ a CA is simply denoted as $\CA(N;t,k,v)$ \cite{Handbook}.

CAs are known under this name since a study of Sloane \cite{Sloane}, according to whom they have been first studied by R\'enyi for strength $t=2$ and alphabet size $v=2$ in \cite{renyi2007foundations}.
Covering arrays have been studied under different names, such as surjective matrices \cite{Honkala1992} or surjective arrays \cite{GreedyCodes}.
Further, covering arrays over binary alphabets are equivalent to independent families of sets \cite{KleitmanSpencer}, \cite{binarySurvey},
and transversal covers as considered in \cite{stevens_lower_1998} are covering arrays of strength two.
However, CAs are still subject to current research, which includes exploration using computational methods \cite{kokkalaStrength2structure},
connections to related combinatorial designs such as \emph{covering perfect hash-families} \cite{colbourn_asymptotic_2018},
and generalizations to \emph{variable strength covering arrays} \cite{LuciaVCA}.
A recent survey of CAs and their constructions can be found in \cite{PepeSurvey}.
The combinatorial aspects of CAs are discussed in \cite{Aspects}, where several constructions and relations to other combinatorial designs, such as difference matrices are presented.

\paragraph{Applications of Covering Arrays}
In recent years CAs have attracted much attention from practitioners due to their application in a branch of software testing, referred to as \emph{combinatorial testing} \cite{kuhn2013introduction},
where the rows of CAs give rise to software tests \cite{HartmanSoftware, HartmanProblems}.
The characteristic of CAs that every $t$-tuple occurs in every $N\times t$ subarray as a row has shown to be especially beneficial for software testing. With this property CAs can be used to reveal software failures resulting from $t$-way interactions of input parameters \cite{InputSpaceCoverageMatters}.
Empirical studies suggest that in some cases the application of combinatorial testing based on CAs of strength $t=4$ to $t=6$ approaches pseudo-exhaustive testing. More specifically, failure report analyses \cite{kuhn_estimating_2016} reveal that most failures are caused by combinations of at most four to six input parameters, depending on the use case.
However, applications of CAs are not merely limited to software testing, as
they are further applicable for hardware testing \cite{HTpaper}, network testing \cite{558835}, and with enhanced properties (as \emph{locating arrays}), they can be used for screening of wireless networks \cite{LAforScreeningComplexEngineeredSystems}.
All these applications have in common, that the deployed tests are based on the rows of a CA.
%
Since practitioners are interested in minimizing resources used for testing, they are interested in using CAs with a \emph{small} number of rows for their testing processes.

This is where the interests of practitioners and theorists meet: in the generation of CAs with as few rows as possible.
The challenge of this problem becomes apparent when considering the plurality of methods, constructions and algorithms,
that are required to constitute the state of the art of smallest CAs, i.e. $\CA(N;t,k,v)$ with \emph{minimized} $N$ for given $t,k,v$.
An overview of the state of the art is maintained in the covering array tables by Colbourn \cite{CJtables} (currently unavailable),  accessible under \cite{CJtablesWayBack}.

\paragraph{Optimal Covering Arrays and Covering Array Numbers}
The covering array number $\CAN_{\lambda}(t,k,v)$ is defined as the minimal number $N$, such that a $\CA_{\lambda}(N;t,k,v)$ exists.
Again, for $\lambda = 1$ we denote the covering array number (\emph{CAN} for short) as $\CAN(t,k,v)$.
A CA with a minimal number of rows is called \emph{optimal}.
Similarly to CAN, for given $N,t,v$ the maximal number of columns $k$ where a $\CA(N;t,k,v)$ exists is denoted as $\CAK(N;t,v)$.
These two numbers are clearly connected via the relations $\CAN(t,k,v)=\min\{N : \CAK(N;t,v) \geq k \}$ and
$\CAK(N;t,v)=\max\{k : \CAN(t,k,v) \leq N \}$, see also \cite{newCANs} or \cite{binarySurvey}.
It is worth to mention that, to the best of the authors' knowledge, all known CANs have been established in a constructive manner,
in the sense that a construction or computation of an optimal CA was used to establish them.
In this sense, as of this writing it is tantamount to determine an \emph{optimal CA} and to determine a \emph{CAN}. 

In general, the problem of finding optimal CAs remains unsolved, however, there are some special cases where they are known.
As such, optimal CAs $\CA(N;2,k,2)$  are known for all $k \in \mathbb N$ due to a connection to maximal anti-chains and sperner-type theorems,
examined by Katona \cite{KatonaSpernerTypeTheorems} as well as Kleitman and Spencer \cite{KleitmanSpencer}.
%
Therefore also the value of $CAN(2,k,2)$ is known for $k\geq 2$ to be
\begin{equation}
	\label{eq:binaryStr2}
	\CAN(2,k,2) = \min_{N \in \mathbb N} \left\{k \leq \binom{N-1}{\left\lceil \frac{N}{2} \right\rceil} \right\}.
\end{equation}

Further, some optimal CAs are known as they coincide with orthogonal arrays of index one (in general we have $\OA(\lambda v^t,k,v,t) = \CA_{\lambda}(\lambda v^t;t,k,v)$), which can for example be constructed using Bush's construction \cite{bush1952orthogonal-habil}.
An introduction to orthogonal arrays, including information on their existence and their properties, is given in \cite{SloaneOA}, where we refer the interested reader for more information.
Aside from these theoretical constructions that are based on connections to other structures appearing in discrete mathematics, which typically yield a family of optimal CAs,
only a few other optimal CAs are known.
%
Some optimal CAs have been determined based on others by making use of bounds on CAN,
such as the elementary bound
\begin{equation}
	\label{eq:bound}
	\CAN(t-1,k-1,v)\leq \frac{1}{v}\CAN(t,k,v),
\end{equation}
which is also mentioned in \cite{Aspects}, and is frequently used in other works to establish new CAN values, or to strengthen bounds on CAN as in \cite{CJclassification}.
However, new results on optimal CAs are mostly established with the aid of computational methods, which we will review in the next section.
%

\paragraph{Contribution}
In this article, the method of combining theoretical properties of \emph{balanced covering arrays}, revisited in Section \ref{sec:balCAs}, with a result from computational classification is used for the first time to establish new CANs.
Using this method, we show the lower bounds $18 \leq \CAN(3,k,2)$ for $k \in {17,18,19,20}$ and $36 \leq \CAN(4,18,2)$,
which determine the respective CANs and establish optimality of five (previously known) covering arrays.

\section{Related Work}

Classification of CAs, i.e. enumeration of non-equivalent CAs, is a useful technique to establish non-existence results.
Known equivalence actions of CAs, i.e. actions that preserve their defining properties, are permutations of columns,
permutations of rows and permutations of the symbols of a column. 
See \cite{balCApaper} for more details.
Classification of CAs has been used in several of the following works related to optimal CAs.

Upper bounds $\CAN(t,k,v)\leq N$ are usually shown by generating a \allowbreak $\CA(N;t,k,v)$, while lower bounds on CANs are often established computationally by complete or exhaustive search algorithms.
%
An overview of upper and lower bounds on CAN for CAs with up to $k=14$ columns is given by Colbourn et al. \cite{CJclassification}. 
The authors derived these bounds using various techniques, ranging from theoretical results and combinatorial constructions to exhaustive search and metaheuristic search techniques.
Additionally, the authors present classification results whenever feasible.
However, since their work \cite{CJclassification} was published some additional CANs and bounds on CANs have been determined.

Chateauneuf,  Colbourn and Kreher \cite{CAsofStrengthThree} show $\CAN(3,6,3)=33$ by using a combination of a computational result, a theoretical bound and a combinatorial construction, as follows.
The bound $\CAN(2,5,3)\geq 11$ is a computational result determined with linear integer programming that was reported in \cite{Sloane}.
Combining this with inequality \eqref{eq:bound} yields $33\leq\CAN(3,6,3)$.
Further, the mathematical construction presented in \cite[Theorem 2.1]{CAsofStrengthThree} allows to derive a $\CA(33;3,6,3)$, which establishes $\CAN(3,6,3)=33$.

The work of Hnich et al. \cite{HnichCP} demonstrates how CAN values or bounds thereof can be determined computationally.
They present a constraint satisfaction problem (CSP) formulation for CA generation, which, together with a complete CSP solver, allows determining the existence of a CA for given parameters $N,t,k$ and $v$.
Using the presented CSP formulation together with a complete CSP solver, optimality of nine CAs was shown,
while for several other parameters $t, k$ and $v$ the lower or upper bound on CAN was strengthened.

Izquierdo-Marquez and Torres-Jimenez \cite{constrNonIsoCAs} present a classification algorithm where for a given number of rows $N$,
strength $t$ and alphabet size $v$ all non-equivalent CAs are enumerated.
For this purpose, an array is extended with additional columns, in combination with backtracking when the array cannot be further extended.
Using the computational results from this algorithm they present five new CANs.
In an extension of their work the same authors derive additional twelve new CANs with an improved version of their algorithm \cite{orderlyAlgorithm}.
Additionally, based on the computational result $\CAN(2,13,3)=15$ and inequality \eqref{eq:bound} they show optimality of $\CA(45;3,14,3)$, i.e. $\CAN(3,14,3)=45$.

%
The same authors propose further a different type of classification algorithm in \cite{newCANs}, where
they use the construction underlying the proof of the inequality \eqref{eq:bound}, see e.g., the argument for equation (4) in \cite{Aspects}.
In this approach all non-equivalent CAs $\CA(N;t+1,k+1,v)$ are generated by a vertical juxtaposition of different CAs $\CA(N_i,t,k,v)$ of reduced strength, for $i=1,\ldots,v$.
Using this classification algorithm, several lower bounds on CAN could be improved, and some CAN values could be determined.
For example it was shown that $\CAN(4,13,2)=32$.
In combination with inequality \eqref{eq:bound} additionally the CAN values $\CAN(5,14,2)= 64$, $\CAN(6,15,2)= 128$ and $\CAN(7,16,2)= 256$ were established,
since the lower bound gained from inequality \eqref{eq:bound} already met the known upper bound.
The same inequality was applied to compute improved lower bounds for $\CAN(t,t+3,2)$ for $t=7,\ldots,11$ from the computational result $\CAN(6,9,2)\geq 107$.
In total six new CANs and $13$ improved lower bounds on CANs were reported \cite{newCANs}.

Kokkala et al. \cite{kokkalaStrength2structure} present a classification algorithm specifically for strength $t=2$ CAs.
Their algorithm also uses column extension and backtracking, 
together with an additional pruning condition to restrict the search space when searching for CAs with $N$ rows and alphabet size $v$, where $N<v(v+1)$.
Further, double-counting (based on the orbit‐stabilizer theorem) is used to increase the confidence in the computational results.
For finding columns for extending the considered array and for identifying non-equivalent CAs,
the authors represent these problems as graph problems as well as linear equation systems and use respective solvers to derive solutions.
The authors report on determining 10 new CANs.

\subsection{Bounds on $\CAN(3,20,2)$ and $\CAK(17;3,2)$}
For the CAN values established in this paper, $\CAN(3,20,2)$ plays a central role.
We briefly outline a chronology of upper and lower bounds for this value.

Sloane studied the connection between CAs of strength three and intersecting codes, reporting an upper bound of $18$ for $\CAN(3,20,2)$ in \cite[Tab. III]{Sloane} in 1993,
giving reference to a construction due to Roux \cite{RouxThesis}, which is nowadays known together with its generalizations as \emph{Roux-type constructions} \cite{CJroux-type}. Since any $k-1$ columns of a $\CA(N;t,k,v)$ constitute a $\CA(N;t,k-1,v)$ as long as $k-1 \geq t$, this immediately yields upper bounds $\CAN(3,k,2)\leq 18$ for $3\leq k \leq 20$.
In the same work Sloane strengthens these bounds for all $3\leq k \leq 16$ to $\CAN(3,k,2) \leq 17$ by explicitly constructing a $\CA(17;3,16,2)$ based on a $2 - (16,6,2)$ biplane, attributing the construction to  D. L. Kreher and V. D. Tonchev.
We like to mention that this does not impose any lower bound on $\CAN(3,k,2)$ for $k \in \{17,18,19,20\}$, since the upper bound on $\CAN(3,k,2)$ for $k\leq 16$ does not exclude the existence, e.g., of a $\CA(16;3,20,2)$.

Chateauneuf and Kreher \cite[Tab. V]{ChateauneufKreher} additionally report the lower bound $16 \leq \CAN(3,k,2) \leq 18$ for $k \in \{17,18,19,20\}$ in their work published in 2002. This lower bound can be derived from $8 \leq \CAN(2,k-1,2)$ for $k \in \{17,18,19,20\}$ using \eqref{eq:bound}; we refer to \cite{ChateauneufKreher} for more details.

In 2006 Colbourn et al. \cite{CJroux-type} devised various \emph{Roux-type} constructions for CAs, particularly of strength three and four. Their work includes a rigorous treatment of the numerical consequences for bounds on CANs due to the established constructions.
In the plethora of reported upper bounds on CAN, we can find again the upper bound $\CAN(3,k,2)\leq 18$ for $k \in \{17,18,19,20\}$, see \cite[Tab. 1]{CJroux-type}.

More recent works, such as Colbourn et al. 2010 \cite{CJclassification}, as well as Izquierdo-Marquez and Torres-Jimenez 2018 \cite{orderlyAlgorithm} and 2019 \cite{newCANs} do not treat the case of interest, seemingly due to computational limitations.
For example, the computational results of \cite{orderlyAlgorithm} show the non-existence of any $\CA(16;3,15,2)$ and therefore optimality of any $\CA(17;3,16,2)$.
This result strengthens the lower bound on $\CAN(3,k,2)$ to $17 \leq \CAN(3,k,2)$ for $k \in \{ 17,18,19,20 \}$.

To the best of the authors' knowledge, the bounds $17 \leq \CAN(3,k,2) \leq 18$ for $k \in \{17,18,19,20\}$,
are the tightest known as of this writing.

When we switch the formulation to the dual \textit{CAK} problem, we ask how many columns in a binary CA of strength three can be attained with a fixed number of rows.
Thus $\CAN(3,k,2) \leq 17$ for $3\leq k \leq 16$ tells us that $\CAK(17;3,2) \geq 16$,
and $\CAN(3,k,2)\leq 18$ for $3\leq k \leq 20$ tells us that $\CAK(18;3,2) \geq 20$.
%
%
However, we do not know whether $\CAK(17;3,2) > 16$. 

\section{Revisiting Balanced Covering Arrays}\label{sec:balCAs}

In \cite{balCApaper} \emph{balanced CAs} were introduced as intersections of classes of covering arrays (CAs) and packing arrays (PAs).
Such balanced CAs obey to upper and lower bounds regarding the appearance of tuples in subarrays.
More detailed, balanced CAs are defined as follows:

\begin{definition}[Balanced CA, \cite{balCApaper}]
A $\CA(N;t,k,v)$ is called $(\balpha,\bbeta)$-balanced with $\balpha=(\alpha_1,\ldots,\alpha_t) \in \mathbb N^t$ and $\bbeta=(\beta_1,\ldots,\beta_t)  \in \mathbb N^t$ if and only if
for all $i \in \{1,\ldots,t\}$ in each $N \times i$ subarray, each $i$-tuple over the alphabet of size $v$ appears at least $\alpha_i$ and at most $\beta_i$ times as a row of the subarray.
Such an array is denoted as $\CA_{\balpha}^{\bbeta}(N;t,k,v)$.
\end{definition}

Further, the maximal number of columns for which a $(\balpha,\bbeta)$-balanced CA exists is defined as
$\CAK_{\balpha}^{\bbeta}(N;t,v):= \max \{k\in \mathbb N: \exists \CA_{\balpha}^{\bbeta}(N;t,k,v) \}$.

The work in \cite{balCApaper} presents classification results for balanced CAs
which yield values of $\CAK_{\balpha}^{\bbeta}(N;t,v)$ for several parameters $N, t$ and $v$.
These results are obtained by means of an exhaustive search algorithm that makes use of SAT solving or a pseudo‐Boolean
solver.
Given a number of rows $N$ and balance vectors $\balpha$ and $\bbeta$,
the proposed algorithm starts with an empty matrix and uses column extension
with backtracking to generate $(\balpha,\bbeta)$-balanced CAs.
A lex‐leader ordering is used to identify a unique representative of each class of non-equivalent balanced CAs.
To prevent from performing redundant searches, several inequalities or bounds on the balance vectors $\balpha$ and $\bbeta$ are used.
For more details we kindly refer the interested reader to \cite{balCApaper}.
%

%
In the following we repeat some of the inequalities from \cite[Corollary 1]{balCApaper} that are relevant for this work.
In addition to the obvious relation $\lambda_t \geq 1$,
when analyzing balanced CAs it is sufficient to consider balance vectors $\balpha=(\alpha_1,\ldots,\alpha_t)$ and $\bbeta=(\beta_1,\ldots,\beta_t)$ satisfying the bounds given below:
%
\begin{eqnarray}
	\alpha_i &\geq& v\cdot \alpha_{i+1}, \text{ for } 1 \leq i <t, \label{eq:searchTrivBound1} \\
	\beta_{i+1} &\leq& \beta_{i} - (v-1) \alpha_{i+1}, \text{ for } 1 \leq i < t, \label{eq:searchTrivBound3} \\
	\beta_i &\leq& N-(v^i-1)\alpha_i, \text{ for } 1 \leq i \leq t. \label{eq:searchTrivBound6}
\end{eqnarray}

Additionally, for CAs with at least $k$ columns, the bound 
\begin{equation}
	\alpha_i \geq \CAN_{\alpha_t}(t-i,k-i,v), \text{ for } 1 \leq i <t \label{eq:searchCAN}
\end{equation}
can be used, which is in many cases stronger than the bound given by inequality \eqref{eq:searchTrivBound1}. For the proof of these bounds we refer to \cite{balCApaper}.

\section{The non-existence of a $\CA(17;3,17,2)$}
Although in theory the classification algorithm we presented in \cite{balCApaper} is capable of classifying all CAs, in practice this is not possible due to extensive runtimes. For example, we were not able to compute a full classification of balanced CAs $\CA_{\balpha}^{\bbeta}(17;3,k,2)$ for all possible balance vectors $\balpha,\bbeta$. To extend the results published in \cite{balCApaper} we executed an additional computation for $\CA_{_{(8,2,1)}}^{^{(9,7,6)}}(17;3,k,2)$ without a time limit.
This (single-threaded) computation, which is only a part of the one required for computationally classifying all $\CA(17;3,k,2)$, finished only after about $12284$ hours, almost $512$ days.
The results of this computation show that three non-equivalent such balanced CAs with $16$ columns exist and further that there is no $\CA_{_{(8,2,1)}}^{^{(9,7,6)}}(17;3,k,2)$ with $k=17$ columns, therefore $\CAK_{_{(8,2,1)}}^{^{(9,7,6)}}(17;3,2) = 16$.
As we show below,
the inequalities for balanced CAs given in \cite{balCApaper} allow to generalize this non-existence result from balanced CAs with $k=17$ columns to CAs with $k=17$ columns.

\begin{theorem}\label{rem:CAN17}
	$\CAK(17;3,2)=16$.
\end{theorem}
\proof
As mentioned above,
we know $\CAK_{_{(8,2,1)}}^{^{(9,7,6)}}(17;3,2) = 16$ from a computational classification.
For determining $\CAK(17;3,2)$ it is sufficient to determine $\CAK_{\balpha}^{\bbeta}(17;3,2)$ for balance vectors $\balpha$, $\bbeta$ that do not impose any balance constraints on the covering 
arrays.
As such, from the basic inequality in \eqref{eq:searchTrivBound1} we get $\balpha \geq (4,2,1)$ and in combination with \eqref{eq:searchTrivBound6} we get $\bbeta \leq (13,11,10)$,
thus $\CAK(17;3,2) = \CAK_{_{(4,2,1)}}^{^{(13,11,10)}}(17;3,2)$, as these balance vectors do not impose any effective balance constraints.
However, we can further strengthen the upper and lower bounds of the balance vectors without imposing balance constraints.
In particular, when determining the existence of a $\CA_{\balpha}^{\bbeta}(17;3,k,2)$ with $k\geq 17$ inequality \eqref{eq:searchCAN} implies that $\lambda_1 \geq \CAN_{\lambda_t}(2,16,2) \geq 8$, since for $\lambda_t = 1$ the value $\CAN(2,16,2) = 8$ is known due to \cite{KleitmanSpencer} and \cite{KatonaSpernerTypeTheorems}, see also inequality \eqref{eq:binaryStr2}.
Together with inequality \eqref{eq:searchTrivBound6} and $\balpha \geq (4,2,1)$, this yields also the bound $\beta_1\leq N - (v^i-1)\lambda_1 = 9$ for $\beta_1$,
and with \eqref{eq:searchTrivBound3} we obtain the bounds  for $\beta_2$ and for $\beta_3$:
$\beta_2 \leq \beta_1 - (v-1)\lambda_{2} = 7$, and $\beta_3 \leq \beta_2 - (v-1)\lambda_{3} = 6$.
Together we obtain that the balance vectors $\balpha=(8,2,1)$ and $\bbeta=(9,7,6)$ do not impose effective constraints on the balance, and thus on the existence, of a $\CA_{\balpha}^{\bbeta}(17;3,k,2)$ with $k\geq 17$ columns.
This means that any $\CA(17;3,k,2)$ is also a $\CA_{_{(8,2,1)}}^{^{(9,7,6)}}(17;3,k,2)$, for $k\geq 17$,
and further
\begin{eqnarray}
	\text{for } k\geq 17: \CAK(17;3,2)= k \Leftrightarrow \CAK_{_{(8,2,1)}}^{^{(9,7,6)}}(17;3,2) = k.
	\label{eq:helperCAN17}
\end{eqnarray}
From the computational search we get $\CAK_{_{(8,2,1)}}^{^{(9,7,6)}}(17;3,2)=16$,
and hence from the equivalence in \eqref{eq:helperCAN17} we know that $\CAK(17;3,2)$ cannot be $17$ or larger.
It follows that $\CAK(17;3,2) = 16$.

\section{The Optimality of $\CA(18;3,20,2)$ and $\CA(36;4,18,2)$}
\begin{corollary}
	$\CAN(3,k,2)=18$ for $k \in \{17,18,19,20\}$.
\end{corollary}
\proof
Due to $\CAN(t,k,v)=\min\{N : \CAK(N;t,v) \geq k \}$ (see also \cite{newCANs}), a consequence of Theorem \ref{rem:CAN17} is
$\CAN(3,17,2) = \min\{N : \CAK(N;3,2) \geq 17 \} > 17 $.
Further, it is known that $\CAN(3,20,2) \leq 18$, see Sloane \cite[Tab. III]{Sloane}, respectively Theorem 6 therein due to Roux \cite{RouxThesis}.
Consequentially, we have $17 < \CAN(3,17,2) \leq \CAN(3,18,2) \leq \CAN(3,19,2) \leq \CAN(3,20,2) \leq 18$,
which shows the claim. \qed

An optimal $\CA(18;3,20,2)$ is given in Figure \ref{fig:optimalCA18}.

\begin{figure}
	\begin{eqnarray*}
	\footnotesize
		\begin{array}{cccccccccccccccccccc}
			0&0&0&0&0&0&0&0&0&0&0&0&0&0&0&0&0&0&0&0\\
			0&0&0&0&0&1&1&1&1&1&0&0&0&0&0&1&1&1&1&1\\
			0&0&1&1&1&0&0&0&1&1&0&0&1&1&1&0&0&0&1&1\\
			0&1&0&1&1&0&1&1&0&0&0&1&0&1&1&0&1&1&0&0\\
			0&1&1&0&1&1&0&1&0&1&0&1&1&0&1&1&0&1&0&1\\
			0&1&1&1&0&1&1&0&1&0&0&1&1&1&0&1&1&0&1&0\\
			1&0&0&1&1&1&1&0&0&1&1&0&0&1&1&1&1&0&0&1\\
			1&0&1&0&1&0&1&1&1&0&1&0&1&0&1&0&1&1&1&0\\
			1&0&1&1&0&1&0&1&0&0&1&0&1&1&0&1&0&1&0&0\\
			1&1&0&0&1&1&0&0&1&0&1&1&0&0&1&1&0&0&1&0\\
			1&1&0&1&0&0&0&1&1&1&1&1&0&1&0&0&0&1&1&1\\
			1&1&1&0&0&0&1&0&0&1&1&1&1&0&0&0&1&0&0&1\\
			0&0&0&0&0&0&0&0&0&0&1&1&1&1&1&1&1&1&1&1\\
			0&0&0&0&1&1&1&1&1&1&1&1&1&1&0&0&0&0&0&0\\
			0&1&1&1&0&0&0&1&1&1&1&0&0&0&1&1&1&0&0&0\\
			1&0&1&1&0&1&1&0&0&1&0&1&0&0&1&0&0&1&1&0\\
			1&1&0&1&1&0&1&0&1&0&0&0&1&0&0&1&0&1&0&1\\
			1&1&1&0&1&1&0&1&0&0&0&0&0&1&0&0&1&0&1&1\\
		\end{array}
	\end{eqnarray*}
	\caption{An optimal $\CA(18;3,20,2)$.}
	\label{fig:optimalCA18}
\end{figure}

\begin{corollary}\label{rem:CAN_t4k18v2}
	$\CAN(4,18,2)=36$.
\end{corollary}
\proof
Inserting the result $\CAN(3,17,2)=18$ into inequality \eqref{eq:bound} we get
	$18=\CAN(3,17,2)\leq \frac{1}{2}\CAN(4,18,2)$,
and therefore $\CAN(4,18,2)\geq 36$. The corresponding upper bound $\CAN(4,18,2)\leq36$ is known due to a CA found using Simulated Annealing by Torres-Jimenez and Rodriguez-Tello \cite{PepeSA}. \qed

\section{Conclusion and Outlook}

In this work we showed how balance properties of CAs can be used to establish lower bounds on CANs which subsequently give optimality results for five CAs.
In particular, using a new computational result in combination with theoretical knowledge about the interplay of balance vectors and the existence of (balanced) covering arrays we established that the covering arrays $\CA(18;3,17,2)$, $\CA(18;3,18,2)$, $\CA(18;3,19,2)$ and $\CA(18;3,20,2)$,
known since Sloane \cite{Sloane}, respectively Roux \cite{RouxThesis}, are optimal.
Thus the entry $N=18$ for $\CA(N;3,k,2)$ with $k\in \{17,18,19,20\}$, in the Covering Array Tables \cite{CJtables}, accessible under \cite{CJtablesWayBack}, constitutes a Covering Array Number. 
Additionally, using an elementary bound on CAN together with the newly established result, we further show optimality of the covering array $\CA(36;4,18,2)$,
which was previously generated by Torres-Jimenez and Rodriguez-Tello using Simulated Annealing \cite{PepeSA}.

We hope that the deployed connections between balance vectors of balanced CAs and (optimal) CAs can be extended in the future,
which may allow for more efficient computations of balanced CAs or to establish new results regarding CA optimality.

%

%

  \bibliographystyle{elsarticle-num} 


\end{document}